\documentclass[12pt,a4paper]{article}
\usepackage{amssymb,amsmath}
\usepackage{amsfonts}
\usepackage{a4}
\begin{document}
 \baselineskip18pt
 \title{\bf  Quadratic residue codes over the ring $\mathbb{F}_{p}[u]/\langle u^m-u\rangle$ and their Gray images}
\author{Mokshi Goyal and Madhu Raka
 \\ \small{\em Centre for Advanced Study in Mathematics}\\
\small{\em Panjab University, Chandigarh-160014, INDIA}\\
\date{}}
\maketitle
 {\abstract{Let $m\geq 2$ be any natural number and let $\mathcal{R}=\mathbb{F}_{p}+u\mathbb{F}_{p}+u^2\mathbb{F}_{p}+\cdots+u^{m-1}\mathbb{F}_{p}$  be a finite non-chain ring, where  $u^m=u$ and $p$ is a prime congruent to $1$ modulo $(m-1)$. In this paper we study  quadratic residue codes over the ring $\mathcal{R}$ and their extensions. A gray map from $\mathcal{R}$ to $\mathbb{F}_{p}^m$  is defined which preserves  self duality of linear codes. As a consequence  self dual, formally self dual and self orthogonal codes are constructed. To illustrate this several examples of   self-dual, self orthogonal and formally self-dual codes are given. Among others a  [9,3,6] linear code over $\mathbb{F}_{7}$ is constructed which is self-orthogonal as well as nearly MDS. The best known linear code with these parameters (ref. Magma) is not self orthogonal.}\vspace{2mm}
 }
 \section{ Introduction.}

Quadratic residue codes are a special kind of cyclic codes of prime length introduced to construct self-dual codes by adding an overall parity check. While initially quadratic residue codes were studied within the confines of finite fields, there have been recent developments on quadratic residue codes over some special rings. Pless and Qian [8] studied quadratic residue codes over $\mathbb{Z}_4$, Chiu et al [3] extended the ideas to the ring $\mathbb{Z}_8$ and Taeri [10] considered QR-codes over $\mathbb{Z}_9$.   Kaya  et al. [4] and Zhang  et al [11] studied quadratic residue codes over $\mathbb{F}_p+u\mathbb{F}_p$.
  Kaya et al [5] studied quadratic residue codes over  $\mathbb{F}_2+u\mathbb{F}_2+u^2\mathbb{F}_2$ whereas Liu et al [6] studied them over non-local ring $\mathbb{F}_p+u\mathbb{F}_p+u^2\mathbb{F}_p$ where $u^3=u$ and $p$ is an odd prime. The authors [9] along with Kathuria   extended their results over the ring
$\mathbb{F}_{p}+u\mathbb{F}_{p}+u^2\mathbb{F}_{p}+u^3\mathbb{F}_{p}$ where $u^4=u$ and $p\equiv 1 ({\rm mod ~} 3)$. Bayram and Siap [1] considered constacyclic codes over the ring $\mathbb{F}_{p}[u]/\langle u^p-u\rangle$. \vspace{2mm}

  This paper is in continuation  of our earlier paper [9]. Here we study, for any natural number $m\geq 2$,  quadratic residue codes and their extensions over the ring $\mathbb{F}_{p}+u\mathbb{F}_{p}+u^2\mathbb{F}_{p}+\cdots+u^{m-1}\mathbb{F}_{p}$ where  $u^m=u$ and $p$ is a prime satisfying $p\equiv 1 ({\rm mod~}(m-1))$, thus  generalizing all the earlier known results.  For any  nonsingular matrix $V$ of order $m$, a Gray map from $\mathcal{R}^n$ to  $\mathbb{F}_p^{mn}$ is defined which preserves linearity.  Further it is shown that if the matrix $V$  satisfies $VV^T=\lambda I_m$ for some $\lambda \in \mathbb{F}_p^*$, where $V^T$ denotes the transpose of the matrix $V$ and $I_m$ denotes the identity matrix of order $m$, then the Gray image  preserves  self duality also. As a consequence  self dual, formally self dual and self orthogonal codes are constructed. To illustrate this several examples of   self-dual, self orthogonal and formally self-dual codes are given such as  a self-dual [36,18,9]  code over $\mathbb{F}_{5}$, a self-dual [12,6,4]  code over $\mathbb{F}_{7}$,  a formally self-dual [36,18,6] code over $\mathbb{F}_{11}$ and   a self orthogonal [33,15,10] code over $\mathbb{F}_5$.   We construct  a self-orthogonal as well as  nearly MDS [9,3,6] code over $\mathbb{F}_{7}$ (nearly MDS in the sense of [2]) which is better than the best known code with these parameters. The dual of our code has parameters [9,6,3]. The [9,6,3] code over $\mathbb{F}_{7}$ constructed by  Liu et al [6] is almost MDS but  not nearly MDS as its dual has parameters [9,3,3]. Further the best known [9,3,6] code over $\mathbb{F}_{7}$ (ref. Magma) is nearly MDS but not self-orthogonal. \vspace{2mm}

The paper is organized as follows: In section 2, we give some  preliminary results and define the Gray map.  In Section 3, we study quadratic residue codes over $\mathcal{R}$, their extensions and give some of their properties. In Section 4, we give several examples of  self-dual, self-orthogonal and formally self-dual codes over $\mathbb{F}_p$.
\section{Preliminary results}
 Throughout the paper, $\mathcal{R}$ denotes the commutative ring $\mathbb{F}_{p}+u\mathbb{F}_{p}+u^2\mathbb{F}_{p}+\cdots +u^{m-1}\mathbb{F}_{p}$, where $u^m=u$, $m\geq 2$ is a natural number
 and $p$ is a prime $\equiv 1({\rm mod}~(m-1))$. $\mathcal{R}$ is a ring of size $p^m$ and characteristic $p$. For   a primitive element $\alpha$ of $\mathbb{F}_{p}$, take $\xi=\alpha^{\frac{p-1}{m-1}}$, so that $\xi^{m-1}=1, \xi\neq 1$ and $\xi^{m-2}+\xi^{m-3}+\cdots+\xi+1=0$. Let $\eta_1, \eta_2, \eta_3,\cdots,\eta_m$ denote the following elements of $\mathcal{R}$.
 \begin{equation} \begin{array}{ll}\eta_1=1-u^{m-1},\vspace{2mm}\\
\eta_2= ({m-1})^{-1}(u+u^2+\cdots+u^{m-2} +u^{m-1}),\vspace{2mm}\\
\eta_3= ({m-1})^{-1}(\xi u+\xi^2u^2+\cdots+\xi^{m-2}u^{m-2} +u^{m-1}),\vspace{2mm}\\
\eta_4= ({m-1})^{-1}(\xi^2 u+(\xi^2)^2u^2+\cdots+(\xi^2)^{m-2}u^{m-2} +u^{m-1}),\vspace{2mm}\\
\cdots ~~~~~~~~~~~~~\cdots ~~~~~~~~~~~~\cdots ~~~~~~~~~~~~~~\cdots\vspace{2mm}\\
\eta_{m}=({m-1})^{-1}(\xi^{m-2} u+(\xi^{m-2})^2u^2+\cdots+(\xi^{m-2})^{m-2}u^{m-2} +u^{m-1}).
\end{array}\end{equation}
A simple calculation shows that \begin{equation}\eta_i^2=\eta_i, ~\eta_i\eta_j=0~{\rm ~for~} 1\leq i, j \leq m, ~i\neq j~{\rm~ and ~}\sum_{i=1}^m \eta_i=1.\end{equation}
\noindent The decomposition theorem of ring theory tells us that  $\mathcal{R}=\eta_1\mathcal{R}\oplus\eta_2\mathcal{R}\oplus \cdots \oplus\eta_m\mathcal{R}$.\vspace{2mm}

\noindent For a linear code $\mathcal{C }$ of length $n$ over the ring $\mathcal{R}$, let \vspace{2mm}\\
$\mathcal{C }_1= \{ x_1\in \mathbb{F}_{p}^n : \exists ~x_2, x_3, \cdots, x_m \in \mathbb{F}_{p}^n {\rm ~such ~that~} \eta_1x_1+\eta_2x_2+\cdots+\eta_m x_m \in \mathcal{C }\},$\vspace{2mm}\\
$\mathcal{C }_2= \{ x_2\in \mathbb{F}_{p}^n : \exists ~x_1, x_3, \cdots, x_m \in \mathbb{F}_{p}^n {\rm ~such ~that~} \eta_1x_1+\eta_2x_2+\cdots+\eta_mx_m \in \mathcal{C }\}$,\vspace{2mm}\\
$~~~~~~~~\cdots ~~~~~~~~~~~~~\cdots~~~~~~~~~~~\cdots~~~~~~~~~~~~~~\cdots$\vspace{2mm}\\
$\mathcal{C }_m= \{ x_m\in \mathbb{F}_{p}^n : \exists ~x_1,x_2, \cdots, x_{m-1} \in \mathbb{F}_{p}^n {\rm ~such ~that~} \eta_1x_1+\eta_2x_2+\cdots+\eta_mx_m \in \mathcal{C }\}$.\vspace{2mm}\\
Then $\mathcal{C}_i, ~i=1,2,\cdots,m$ are linear codes of length $n$ over $\mathbb{F}_{p}$,  $\mathcal{C}=\eta_1\mathcal{C}_1\oplus\eta_2\mathcal{C}_2\oplus\cdots \oplus\eta_m\mathcal{C}_m
$ and $|\mathcal{C }|= |\mathcal{C }_1|~|\mathcal{C }_2|~\cdots~|\mathcal{C }_m|$.
For a code $\mathcal{C}$  over $\mathcal{R}$, the dual code  $\mathcal{C}^\bot$ is defined as $ \mathcal{C}^\bot =\{x\in \mathcal{R}^n | x\cdot y=0 ~ {\rm for ~ all~} y \in \mathcal{C}\}$ where $x\cdot y$ denotes the usual Euclidean inner product. $\mathcal{C}$ is self-dual if $\mathcal{C}=\mathcal{C}^\bot$ and self-orthogonal if $\mathcal{C}\subseteq \mathcal{C}^\bot$.\\
The following result is a simple generalization of a result of [6].\vspace{2mm}

\noindent{\bf Theorem 1}: Let $\mathcal{C}=\eta_1\mathcal{C}_1\oplus\eta_2 \mathcal{C}_2\oplus \cdots\oplus \eta_m\mathcal{C}_m$ be a linear code of length $n$ over $\mathcal{R}$. Then \vspace{2mm}

\noindent (i)~~ $\mathcal{C}$  is cyclic over $\mathcal{R}$ if and only if $\mathcal{C}_i, ~i=1,2,\cdots,m$ are cyclic over $\mathbb{F}_p$.\vspace{2mm}

\noindent (ii) ~~If $\mathcal{C}_i=\langle  g_i (x)\rangle, ~g_i(x)\in \frac{\mathbb{F}_p[x]}{\langle x^{n}-1\rangle}$, $g_i(x)|(x^n-1)$,  then $\mathcal{C}=\langle \eta_1g_1(x),\eta_2g_2(x),\cdots,\eta_mg_m(x)\rangle$\\ $~~~~~~~~ =\langle g(x)\rangle $ where $g(x)= \eta_1g_1+\eta_2g_2+\cdots+\eta_mg_m$ and $g(x)|(x^{n}-1)$.\vspace{2mm}

\noindent (iii)~~ Further $|\mathcal{C }|=p^{mn-\sum_{i=1}^{m}deg(g_i)}$.\vspace{2mm}

\noindent (v)~~ Suppose that $g_i(x)h_i(x)=x^n-1,~ 1\leq i\leq m.$ Let $ h(x)=\eta_1h_1(x)+\eta_2h_2(x)+$ \\$~~~~~~~~~\cdots+\eta_mh_m(x),$ then
     $g(x)h(x)=x^n-1$. \vspace{2mm}

\noindent(vi)~~ $ \mathcal{C}^\perp=\eta_1\mathcal{C}_1^\perp\oplus\eta_2\mathcal{C}_2^\perp\oplus\cdots
     \oplus\eta_m\mathcal{C}_m^\perp.$ \vspace{2mm}

\noindent (vii)~~ $ \mathcal{C}^\perp=\langle h^\perp(x)\rangle,$
    where $ h^\perp(x)=\eta_1h_1^\perp(x)+\eta_2h_2^\perp(x)+\cdots+\eta_mh_m^\perp(x)$,
      $ h_i^\perp(x)$ \\ $~~~~~~~~~$is the reciprocal
    polynomial of $h_i(x), ~1\leq i\leq m.$  \vspace{2mm}

 \noindent(viii)$~~ |\mathcal{C}^\perp|=p^{\sum_{i=1}^m deg(g_i)}$.\vspace{2mm}

\noindent The following is a well known result :\vspace{2mm}

\noindent{\bf Lemma 1}: (i) Let $C$  be a cyclic code of length $n$ over a finite ring $S$ generated by the idempotent $E$ in $S[x]/\langle x^n-1\rangle$ then $C^{\perp}$ is generated by the idempotent $1-E(x^{-1})$.\vspace{2mm}

\noindent (ii) Let $C$ and $D$ be cyclic codes of length $n$ over a finite ring $S$ generated by the idempotents $E_1, E_2$ in $S[x]/\langle x^n-1\rangle$ then $C\cap D$ and $C+D$ are generated by the idempotents $E_1E_2$ and $E_1+E_2-E_1E_2$ respectively.\vspace{2mm}

\noindent We define a Gray map $\Phi : \mathcal{R}\rightarrow \mathbb{F}_p^m$ ~~ given by $$\begin{array}{ll}r(u)=a_0+a_1u+a_2u^2+\cdots +a_{m-1}u^{m-1}\longmapsto(r(0),r(1),r(\xi),\cdots, r(\xi^{m-2}))V\vspace{4mm}\\
= (a_0,a_1,a_2,\cdots ,a_{m-1})\left(
 \begin{array}{cccccc}
                                                                               1 & 1 & 1& 1& \cdots & 1 \\
                                                                               0 & 1 & \xi & \xi^2 &\cdots & \xi^{m-2} \\
                                                                                0 & 1 & \xi^2 & (\xi^2)^2 &\cdots & (\xi^2)^{m-2} \\
                                                                                0 & 1 & \xi^3 & (\xi^3)^2 &\cdots & (\xi^3)^{m-2} \\
                                                                                \cdots & \cdots&\cdots&\cdots&\cdots&\cdots\\
                                                                                 0 & 1 & \xi^{m-2} & (\xi^{m-2})^2 &\cdots & (\xi^{m-2})^{m-2} \\
                                                                                 0&1&1&1&\cdots&1
                                                                             \end{array}
                                                                           \right)V\vspace{4mm}\\
                                                                            =(a_0,a_1,a_2,\cdots ,a_{m-1})M V\end{array}
$$ where $M$ is an $m\times m$ nonsingular matrix of Vandermonde determinant $\hspace{-2mm}{\displaystyle\prod_{1\leq i<j\leq m-1}}\hspace{-4mm}{(\xi^j-\xi^i)}$ and  $V$ is any nonsingular matrix over $\mathbb{F}_p$ of order $m\times m$.
 This map can be extended from $\mathcal{R}^n$ to  $\mathbb{F}_p^{mn}$ component wise.\vspace{2mm}

 \noindent Define the Gray weight of a an element $r \in \mathcal{R}$ by $w_{G}(r) =w_H(\Phi(r))$, the Hamming weight of $\Phi(r)$. The Gray weight of  a codeword
$c=(c_0,c_1,\cdots,c_{n-1})\in \mathcal{R}^n$ is defined as $w_{G}(c)=\sum_{i=0}^{n-1}w_{G}(c_i)=\sum_{i=0}^{n-1}w_H(\Phi(c_i))=w_H(\Phi(c))$. For any two elements $c_1, c_2 \in \mathcal{R}^n$, the Gray distance $d_{G}$ is given by $d_{G}(c_1,c_2)=w_{G}(c_1-c_2)=w_H(\Phi(c_1)-\Phi(c_2))$. \vspace{2mm}

 \noindent \textbf{Theorem 2.} The Gray map $\Phi$ is an  $\mathbb{F}_p$- linear, one to one and onto map. It is also distance preserving map from ($\mathcal{R}^n$, Gray distance $d_{G}$) to ($\mathbb{F}_p^{mn}$, Hamming distance). Further if the matrix $V$  satisfies $VV^T=\lambda I_m$, $\lambda \in \mathbb{F}_p^*$, where $V^T$ denotes the transpose of the matrix $V$, then the Gray image $\Phi(\mathcal{C})$ of a self-dual code $\mathcal{C}$  over $\mathcal{R}$ is a self-dual code in $\mathbb{F}_p^{mn}$. \vspace{2mm}

\noindent \textbf{Proof.} The first two assertions holds as $MV$ is an invertible matrix over $\mathbb{F}_p$.\\ Let now $V=(v_{ij})$, $1\leq i,j\leq m$, satisfying $VV^T=\lambda I_m$. So that
 \begin{equation} {\displaystyle \sum_{k=1}^m}~ v_{jk}^2=\lambda ~~{\rm for ~all ~} j, 1\leq j \leq m {\rm ~~and ~~}  {\displaystyle \sum_{k=1}^m}~ v_{jk}v_{\ell k}=0 ~~{\rm for ~} j \neq \ell. \end{equation}
 Let $\mathcal{C}$ be a self dual code over $\mathcal{R}$. Let $r=(r_0,r_1,\cdots,r_{n-1}), s=(s_0,s_1,\cdots,s_{n-1}) \in \mathcal{C}$ where $r_i=a_{i0}+a_{i1}u+\cdots +a_{im-1}u^{m-1}$ and $s_i=b_{i0}+b_{i1}u+\cdots +b_{im-1}u^{m-1}$. Then $$0=r\cdot s=\sum_{i=0}^{n-1}r_is_i= \sum_{i=0}^{n-1}~\sum_{j=0}^{m-1}~\sum_{\ell=0}^{m-1}~a_{ij}~b_{i\ell}~u^{j+\ell}$$ implies that (comparing the coefficients of $u^r$
on both sides)\begin{equation}  {\displaystyle \sum_{i=0}^{n-1}}a_{i0}b_{i0}=0 ~~~~~~~~{\rm and}\vspace{-4mm}\end{equation}  \begin{equation}{\displaystyle \sum_{i=0}^{n-1}}(a_{i0}b_{ir}+a_{i1}b_{i,r-1}+\cdots +a_{ir}b_{i0}+a_{ir}b_{i,m-1}+a_{i,r+1}b_{i,m-2}+\cdots+a_{i,m-1}b_{i0})=0,\end{equation} for each $r$, $1\leq r \leq m-1.$\vspace{2mm}

\noindent For convenience we call $(r_{i}(0),r_i(1),r_i(\xi),\cdots,r_i(\xi^{m-2}))=( \alpha_{i1},\alpha_{i2},\cdots,\alpha_{im})$ and
$(s_{i}(0),s_i(1),s_i(\xi),\cdots,s_i(\xi^{m-2}))=( \beta_{i1},\beta_{i2},\cdots,\beta_{im})$. Then
$$\Phi(r_i)=( \alpha_{i1},\alpha_{i2},\cdots,\alpha_{im})V =\big(\sum_{j=1}^{m}~\alpha_{ij}~v_{j1},\sum_{j=1}^{m}~\alpha_{ij}~v_{j2}~,\cdots, \sum_{j=1}^{m}~\alpha_{ij}~v_{jm}\big)$$ Similarly $$\Phi(s_i) =\big(\sum_{\ell=1}^{m}~\beta_{i\ell}~v_{\ell 1},\sum_{\ell=1}^{m}~\beta_{i\ell}~v_{\ell 2}~,\cdots, \sum_{\ell=1}^{m}~\beta_{i\ell}~v_{\ell m}\big).$$  Using (2), we find that

$$\begin{array}{ll}\Phi(r_i)\cdot \Phi(s_i)& ={\displaystyle\sum_{k=1}^m}~{\displaystyle\sum_{j=1}^{m}}~~{\displaystyle\sum_{\ell=1}^{m}}\alpha_{ij}~\beta_{i\ell}~v_{jk}~v_{\ell k}\\&
={\displaystyle\sum_{j=1, \ell=j}^{m}}\alpha_{ij}~\beta_{ij}\Big({\displaystyle\sum_{k=1}^m}~v_{jk}^2\Big)+{\displaystyle\sum_{j=1}^{m}}
{\displaystyle\sum_{\ell=1,\ell\neq j}^{m}}\alpha_{ij}~\beta_{i\ell}\Big({\displaystyle\sum_{k=1}^m}~v_{jk}~v_{\ell k}\Big)\\&=\lambda {\displaystyle\sum_{j=1}^{m}}\alpha_{ij}~\beta_{ij}.
\end{array}$$
Now
$$\begin{array}{ll}\Phi(r)\cdot \Phi(s)&= {\displaystyle\sum_{i=0}^{n-1}}\Phi(r_i)\cdot \Phi(s_i)= \lambda {\displaystyle\sum_{i=0}^{n-1}}~ {\displaystyle\sum_{j=1}^{m}}\alpha_{ij}~\beta_{ij}\\&={\displaystyle\sum_{i=0}^{n-1}}\Big(\alpha_{i1}~\beta_{i1}+{\displaystyle\sum_{j=2}^{m}}~\alpha_{ij}~\beta_{ij}\Big)
\\&={\displaystyle\sum_{i=0}^{n-1}}a_{i0}b_{i0}+{\displaystyle\sum_{i=0}^{n-1}}~{\displaystyle\sum_{j=2}^{m}}r_i(\xi^{j-2})~s_i(\xi^{j-2})\vspace{2mm}\\
&={\displaystyle\sum_{i=0}^{n-1}}~{\displaystyle\sum_{k=0}^{m-2}}r_i(\xi^{k})~s_i(\xi^{k})\vspace{2mm}\\&={\displaystyle\sum_{i=0}^{n-1}}~{\displaystyle\sum_{k=0}^{m-2}}
\Big({\displaystyle\sum_{j=0}^{m-1}}a_{ik}~\xi^{kj}  \Big)~\Big({\displaystyle\sum_{\ell=0}^{m-1}}b_{i\ell}~\xi^{k\ell}  \Big)\vspace{2mm}\\&={\displaystyle\sum_{i=0}^{n-1}}~{\displaystyle\sum_{k=0}^{m-2}}~{\displaystyle\sum_{j=0}^{m-1}}~
{\displaystyle\sum_{\ell=0}^{m-1}}a_{ik}~b_{i\ell}
~\xi^{k(j+\ell)}\vspace{2mm}\\&=A_0+A_1\xi+A_2\xi^2+\cdots+A_{m-2}\xi^{m-2} ~~~~~~{\rm say}.
\end{array}$$
Using (3) and (4), one can check that each $A_i$ is zero, which proves the result.

\section{Main results}

In this section, quadratic residue codes over $\mathcal{R}$ are defined in terms of their idempotent generators. Let $n=q$ be an odd prime such that $p$ is a quadratic residue modulo $q$. Let $Q_q$ and $N_q$ be the sets of quadratic residues and non-residues modulo $q$ respectively. Let
$$ r(x)= {\displaystyle \prod_{r\in Q_q}}(x-\alpha^r), ~~ n(x)= {\displaystyle \prod_{n\in N_q}}(x-\alpha^n)$$
where $\alpha $ is a primitive $q$th root of unity in some extension field of $\mathbb{F}_p$. Following classical notation of [9], let $\mathbb{Q}$, $\mathbb{N}$ be the QR codes generated by $r(x)$ and $n(x)$ over $\mathbb{F}_p$ and $\tilde{\mathbb{Q}}$, $\tilde{\mathbb{N}}$ be the expurgated QR codes generated by $(x-1)r(x)$ and $(x-1)n(x)$ respectively. we use the notation \vspace{2mm}\\
$~~~~~~j_1(x)=\sum_{i\in Q_q}x^i , ~ j_2(x)=\sum_{i\in N_q}x^i,$\vspace{2mm}\\ $~~~~~~h(x)=1+j_1(x)+j_2(x)=1+x+x^2+\cdots+x^{q-1}=r(x)n(x). $

\noindent If $p>2$,  $q\equiv \pm 1({\rm mod}~4)$ and $p$ is a quadratic residue modulo $q$, then idempotent generators of  $\mathbb{Q}$, $\mathbb{N}$, $\tilde{\mathbb{Q}}$, $\tilde{\mathbb{N}}$ over $\mathbb{F}_p$ are given by, (see [7]),\vspace{2mm}\\
$E_q(x)= \frac{1}{2}(1+\frac{1}{q})+\frac{1}{2}(\frac{1}{q}-\frac{1}{\theta})j_1+\frac{1}{2}(\frac{1}{q}+\frac{1}{\theta})j_2,$\vspace{2mm}\\
$E_n(x)= \frac{1}{2}(1+\frac{1}{q})+\frac{1}{2}(\frac{1}{q}-\frac{1}{\theta})j_2+\frac{1}{2}(\frac{1}{q}+\frac{1}{\theta})j_1,$\vspace{2mm}\\
$F_q(x)= \frac{1}{2}(1-\frac{1}{q})-\frac{1}{2}(\frac{1}{q}+\frac{1}{\theta})j_1-\frac{1}{2}(\frac{1}{q}-\frac{1}{\theta})j_2,$\vspace{2mm}\\
$F_n(x)= \frac{1}{2}(1-\frac{1}{q})-\frac{1}{2}(\frac{1}{q}+\frac{1}{\theta})j_2-\frac{1}{2}(\frac{1}{q}-\frac{1}{\theta})j_1,$\vspace{2mm}\\
respectively where $\theta$ denotes the Gaussian sum and $\chi(i)$ denotes the Legendre symbol, that is
$$ \theta= {\displaystyle \sum_{i=1}^{q-1}}\chi(i)\alpha^i, ~~~~~~\chi(i)= \left\{ \begin{array}{ll}0,& p|i\\1,&i\in Q_q\\-1, & i\in N_q. \end{array}\right.$$
It is known that $\theta^2=-q$ if $q\equiv 3({\rm mod}~4)$ and $\theta^2=q$ if $q\equiv 1({\rm mod}~4).$\\
For convenience we write $d_1=E_q(x),d_2= E_n(x)$, the odd-like idempotents and $ e_1=F_q(x),e_2= F_n(x)$, the even-like idempotents.\vspace{2mm}\\
\noindent{\bf Lemma 2}: $d_1+d_2= 1+\frac{1}{q}h$, $e_1+e_2= 1-\frac{1}{q}h$, $d_1-e_1= \frac{1}{q}h$, $d_2-e_2= \frac{1}{q}h$. Further $d_1d_2= \frac{1}{q}h$ and $e_1e_2=0$.\vspace{2mm}\\
\noindent {\bf Proof} A direct simple calculation gives the first four expressions. To see that $d_1d_2= \frac{1}{q}h$, we note that $d_1d_2= d_1(d_1+d_2-1)=d_1(\frac{1}{q}h)=\frac{1}{q}h$ as $d_1$ is the multiplicative unity of the QR code $\mathbb{Q}$ and $\frac{1}{q}h$ is divisible by $r(x)$ so it belongs to $\mathbb{Q}$. Again $e_1e_2= e_1(e_1+e_2-1)=e_1(\frac{-1}{q}h)=0$ as $e_1 \in
\tilde{\mathbb{Q}}=\langle (x-1)r(x)\rangle$ and $h(x)=r(x)n(x)$, so $e_1h$ is a multiple of $x^q-1$ and hence zero in $\frac{\mathbb{F}_p[x]}{\langle x^{q}-1\rangle}$. \vspace{2mm}

\noindent Let $\mathcal{R}_q$ denote the ring $ \frac{\mathcal{R}[x]}{\langle x^{q}-1\rangle}$. Using the properties (2) of idempotents $\eta_i$, we have \\
\noindent{\bf Lemma 3}: Let $p$ be a prime, $p\equiv 1$(mod $(m-1)$) and $\eta_i, ~1\leq i \leq m $ be idempotents as defined in (1). Then for  $i_1,i_2,\cdots,i_m \in \{1,2\}$ and for any tuple $(d_{i_1}, d_{i_2}, \cdots, d_{i_m})$ of odd-like idempotents not all equal and for any tuple $(e_{i_1}, e_{i_2}, \cdots, e_{i_m})$ of even-like idempotents not all equal,   $\eta_1d_{i_1}+\eta_2d_{i_2}+\cdots+\eta_m d_{i_m}$ and $\eta_1e_{i_1}+\eta_2e_{i_2}+\cdots+\eta_m e_{i_m}$ are respectively odd-like and even-like idempotents in the ring $\mathcal{R}_q= \frac{\mathcal{R}[x]}{\langle x^{q}-1\rangle}$.\vspace{2mm}

We now define quadratic residue codes over the ring $\mathcal{R}$ in terms of their idempotent generators. We denote the set $\{1,2,\cdots,m\}$ by $\mathbb{A}$. For $i \in \mathbb{A}$, let $D_{\{i\}}$ denote the odd-like idempotent of the ring $\mathcal{R}_q$ in which $d_1$ occurs at the $i$th place and $d_2$ occurs at the remaining $ 1,2,\cdots,i-1,i+1,\cdots,m$ places i.e.
\begin{equation}D_{\{i\}}=  \eta_1d_2+\eta_2d_2+\cdots+\eta_{i-1}d_2+\eta_id_1+\eta_{i+1}d_2+\cdots+\eta_md_2= \eta_id_1+(1-\eta_i)d_2.\end{equation}

For $i_1,i_2 \in \mathbb{A}$, $i_1 \neq i_2$ let $D_{\{i_1,i_2\}}$ denote the odd-like idempotent  in which $d_1$ occurs at the $i_1$th and $i_2$th places and $d_2$ occurs at the remaining $ 1,2,\cdots,i_1-1,i_1+1,\cdots,i_2-1,i_2+1,\cdots,m$ places i.e.
\begin{equation}\begin{array}{ll}D_{\{i_1,i_2\}}=&  \eta_1d_2+\eta_2d_2+\cdots+\eta_{i_1-1}d_2+\eta_{i_1}d_1+\eta_{i_1+1}d_2+\cdots+\eta_{i_2-1}d_2+\eta_{i_2}d_1\\&+\eta_{i_2+1}d_2+\cdots+\eta_md_2= (\eta_{i_1}+\eta_{i_2})d_1+(1-\eta_{i_1}-\eta_{i_2})d_2.\end{array}\end{equation}

In the same way, for $i_1,i_2,\cdots, i_k \in \mathbb{A}$,  $i_{r} \neq i_{s},~1\leq r,s \leq k$ let $D_{\{i_1,i_2,\cdots, i_k\}}$ denote the odd-like idempotent
\begin{equation}D_{\{i_1,i_2,\cdots, i_k\}}=(\eta_{i_1}+\eta_{i_2}+\cdots +\eta_{i_k})d_1+(1-\eta_{i_1}-\eta_{i_2}-\cdots -\eta_{i_k})d_2.\end{equation}
For $i \in \mathbb{A}$, $i_1,i_2,\cdots, i_k \in \mathbb{A}$, where $i_{r} \neq i_{s},~1\leq r,s \leq k$ let the corresponding odd-like idempotents be
\begin{equation}D'_{\{i\}}=  \eta_id_2+(1-\eta_i)d_1.\end{equation}
\begin{equation}D'_{\{i_1,i_2,\cdots, i_k\}}=(\eta_{i_1}+\eta_{i_2}+\cdots +\eta_{i_k})d_2+(1-\eta_{i_1}-\eta_{i_1}-\cdots -\eta_{i_k})d_1.\end{equation}
Similarly we define even-like idempotents for $i \in \mathbb{A}$  and $i_1,i_2,\cdots, i_k \in \mathbb{A}$,  $i_{r} \neq i_{s},~1\leq r,s \leq k$,
\begin{equation}E_{\{i\}}= \eta_ie_1+(1-\eta_i)e_2.\end{equation}
\begin{equation}E'_{\{i\}}= \eta_ie_2+(1-\eta_i)e_1.\end{equation}
\begin{equation}E_{\{i_1,i_2,\cdots, i_k\}}=(\eta_{i_1}+\eta_{i_2}+\cdots +\eta_{i_k})e_1+(1-\eta_{i_1}-\eta_{i_2}-\cdots -\eta_{i_k})e_2.~~~~\end{equation}
\begin{equation}E'_{\{i_1,i_2,\cdots, i_k\}}=(\eta_{i_1}+\eta_{i_2}+\cdots +\eta_{i_k})e_2+(1-\eta_{i_1}-\eta_{i_2}-\cdots -\eta_{i_k})e_1.\end{equation}

\noindent Let $Q_{\{i\}},~Q'_{\{i\}},~ Q_{\{i_1,i_2,\cdots, i_k\}},~Q'_{\{i_1,i_2,\cdots, i_k\}}$ denote the odd-like quadratic residue codes  and $S_{\{i\}},~S'_{\{i\}},~ S_{\{i_1,i_2,\cdots, i_k\}},~S'_{\{i_1,i_2,\cdots, i_k\}}$ denote the even-like quadratic residue codes over $\mathcal{R}$ generated by the corresponding  idempotents, i.e. \vspace{2mm}\\
$Q_{\{i\}}= \langle D_{\{i\}}\rangle $,~~
$Q'_{\{i\}}= \langle D'_{\{i\}}\rangle $,~~
$S_{\{i\}}= \langle E_{\{i\}}\rangle $, ~~
$S'_{\{i\}}= \langle E'_{\{i\}}\rangle $,\vspace{2mm}\\
 $Q_{\{i_1,i_2,\cdots, i_k\}}= \langle D_{\{i_1,i_2,\cdots, i_k\}}\rangle $,~~
  $Q'_{\{i_1,i_2,\cdots, i_k\}}= \langle D'_{\{i_1,i_2,\cdots, i_k\}}\rangle $,\vspace{2mm}\\
   $S_{\{i_1,i_2,\cdots, i_k\}}= \langle E_{\{i_1,i_2,\cdots, i_k\}}\rangle $,~~
    $S'_{\{i_1,i_2,\cdots, i_k\}}= \langle E'_{\{i_1,i_2,\cdots, i_k\}}\rangle$.\vspace{2mm}

\noindent{\bf Theorem 3 :} Let $p\equiv 1({\rm mod}~(m-1))$, $q$  an odd prime and $p$ be a quadratic residue modulo $q$. Then for $i \in \mathbb{A}$, $Q_{\{i\}}$ is equivalent to $ Q'_{\{i\}}$ and $S_{\{i\}}$ is equivalent to  $S'_{\{i\}}$. For $i_1,i_2,\cdots, i_k \in \mathbb{A}$,  $i_{r} \neq i_{s},~1\leq r,s \leq k$, $Q_{\{i_1,i_2,\cdots, i_k\}}$ is equivalent to $ Q'_{\{i_1,i_2,\cdots, i_k\}}$, and  $S_{\{i_1,i_2,\cdots, i_k\}}$ is equivalent to $ S'_{\{i_1,i_2,\cdots, i_k\}}$. Further
there are $2^{m-1}-1$ inequivalent odd-like quadratic residues codes and $2^{m-1}-1$ inequivalent even-like quadratic residues codes over the $\mathcal{R}$.\vspace{2mm}

\noindent \textbf{Proof:}
  Let $ n\in N_q$. Let $\mu_n$ be the multiplier map $\mu_n : \mathbb{F}_p \rightarrow \mathbb{F}_p$ given by $\mu_n(a)=an ({\rm mod~} p)$ and acting on polynomials as  $\mu_n(\sum_i f_ix^i)=\sum_if_ix^{\mu_n(i)}$. Then $\mu_n(j_1)=j_2$ and $\mu_n(j_2)=j_1$. Therefore $\mu_n(d_1)=d_2$, $\mu_n(d_2)=d_1$, $\mu_n(e_1)=e_2$, $\mu_n(e_2)=e_1$ and so $\mu_n(\eta_id_1+(1-\eta_i)d_2) =\eta_id_2+(1-\eta_i)d_1$, $\mu_n(\eta_ie_1+(1-\eta_i)e_2) =\eta_ie_2+(1-\eta_i)e_1$, $\mu_n(D_{\{i_1,i_2,\cdots, i_k\}}) =D'_{\{i_1,i_2,\cdots, i_k\}}$, $\mu_n(E_{\{i_1,i_2,\cdots, i_k\}}) =E'_{\{i_1,i_2,\cdots, i_k\}}$. This proves that $Q_i \sim Q'_i$, $S_i \sim S'_i$, $Q_{\{i_1,i_2,\cdots, i_k\}}\sim Q'_{\{i_1,i_2,\cdots, i_k\}}$, and  $S_{\{i_1,i_2,\cdots, i_k\}}\sim S'_{\{i_1,i_2,\cdots, i_k\}}$. \vspace{2mm}

\noindent Note that  $ D_{\mathbb{A}-\{i\}}= D'_{\{i\}}$, $ E_{\mathbb{A}-\{i\}}= E'_{\{i\}}$, $ D_{\mathbb{A}-\{i_1,i_2,\cdots, i_k\}} = D'_{\{i_1,i_2,\cdots, i_k\}}$, $ E_{\mathbb{A}-\{i_1,i_2,\cdots, i_k\}}$ $ = E'_{\{i_1,i_2,\cdots, i_k\}}$. Therefore
\begin{equation} Q_{\mathbb{A}-\{i\}} \sim Q'_{\{i\}} \sim Q_{\{i\}}, ~~ S_{\mathbb{A}-\{i\}}\sim S'_{\{i\}}\sim S_{\{i\}}, \vspace{-2mm}\end{equation} \begin{equation} Q_{\mathbb{A}-\{i_1,i_2,\cdots, i_k\}}\sim Q_{\{i_1,i_2,\cdots, i_k\}},~~ S_{\mathbb{A}-\{i_1,i_2,\cdots, i_k\}} \sim S_{\{i_1,i_2,\cdots, i_k\}}.\end{equation}
 For a given positive integer $k$, the number of choices of the subsets $\{i_1,i_2,\cdots, i_k\}$ of $\mathbb{A}$ is $ m \choose k$.\vspace{2mm}

Let $m$ be even first. Then $|\{i_1,i_2,\cdots, i_{m/2}\}|= |\mathbb{A}-\{i_1,i_2,\cdots, i_{m/2}\}|=\frac{m}{2}$. Using (15) and (16), we find that the number of inequivalent odd-like or even-like QR- codes is $ {m \choose 1} + { m \choose 2}+\cdots { m \choose m/2-1}+ \frac{1}{2}{ m \choose m/2}=2^{m-1}-1$. If $m$ is odd the number of inequivalent odd-like or even-like QR- codes is $ {m \choose 1} + { m \choose 2}+\cdots { m \choose (m-1)/2}=2^{m-1}-1$.\vspace{4mm}

Let $[x]$ denote the greatest integer $\leq x$. we have $[\frac{m}{2}]= \frac{m}{2}$ when $m$ is even and $[\frac{m}{2}]= \frac{m-1}{2}$ when $m$ is odd.\vspace{2mm}

\noindent{\bf Theorem 4 :} If $p\equiv 1({\rm mod}~(m-1))$, $q$ is an odd prime , $p$ is a quadratic residue modulo $q$, then for subsets ${\{i_1,i_2,\cdots, i_k\}}$ of $\mathbb{A} $ with cardinality $k$, $1 \leq k \leq [\frac{m}{2}] $,
the following assertions hold for quadratic residues codes over $\mathcal{R}$.
\vspace{2mm}\\$\begin{array}{ll}

{\rm (i)}& Q_{\{i_1,i_2,\cdots, i_k\}}\cap Q'_{\{i_1,i_2,\cdots, i_k\}}= \langle \frac{1}{q}h\rangle,\vspace{2mm}\\
 {\rm (ii)}& Q_{\{i_1,i_2,\cdots, i_k\}}+Q'_{\{i_1,i_2,\cdots, i_k\}}= \mathcal{R}_q,  \vspace{2mm}\\
 {\rm (iii)}& S_{\{i_1,i_2,\cdots, i_k\}}\cap S'_{\{i_1,i_2,\cdots, i_k\}}= \{0\},\vspace{2mm}\\
 {\rm (iv)}& S_{\{i_1,i_2,\cdots, i_k\}}+ S'_{\{i_1,i_2,\cdots, i_k\}}= \langle 1-\frac{1}{q}h\rangle, \vspace{2mm}\\
{\rm (v)} &S_{\{i_1,i_2,\cdots, i_k\}}\cap \langle \frac{1}{q}h\rangle = \{0\},~S'_{\{i_1,i_2,\cdots, i_k\}}\cap \langle \frac{1}{q}h\rangle = \{0\},\vspace{2mm}\\
 {\rm (vi)} &S_{\{i_1,i_2,\cdots, i_k\}}+ \langle \frac{1}{q}h\rangle = Q_{\{i_1,i_2,\cdots, i_k\}},~
 S'_{\{i_1,i_2,\cdots, i_k\}}+ \langle \frac{1}{q}h\rangle = Q'_{\{i_1,i_2,\cdots, i_k\}},\vspace{2mm}\\
{\rm (vii)}& |Q_{\{i_1,i_2,\cdots, i_k\}}|= p^{\frac{m(q+1)}{2}}, |S_{\{i_1,i_2,\cdots, i_k\}}|= p^{\frac{m(q-1)}{2}}. \end{array}$\vspace{2mm}\\

\noindent \textbf{Proof:}
From the relations (2),(6)-(14) we see that $D_{\{i_1,i_2,\cdots, i_k\}}+D'_{\{i_1,i_2,\cdots, i_k\}}=d_1+d_2 $, $E_{\{i_1,i_2,\cdots, i_k\}}+E'_{\{i_1,i_2,\cdots, i_k\}}=e_1+e_2 $,  $D_{\{i_1,i_2,\cdots, i_k\}}D'_{\{i_1,i_2,\cdots, i_k\}}=d_1d_2 $ and $E_{\{i_1,i_2,\cdots, i_k\}}E'_{\{i_1,i_2,\cdots, i_k\}}=e_1e_2 $.  Therefore by Lemmas 1 and 2, $Q_{\{i_1,i_2,\cdots, i_k\}}\cap Q'_{\{i_1,i_2,\cdots, i_k\}} $ $= \langle D_{\{i_1,i_2,\cdots, i_k\}}D'_{\{i_1,i_2,\cdots, i_k\}}\rangle $ $= \langle \frac{1}{q}h\rangle$,
and  $Q_{\{i_1,i_2,\cdots, i_k\}}+ Q'_{\{i_1,i_2,\cdots, i_k\}}=\langle D_{\{i_1,i_2,\cdots, i_k\}}+D'_{\{i_1,i_2,\cdots, i_k\}}-D_{\{i_1,i_2,\cdots, i_k\}}D'_{\{i_1,i_2,\cdots, i_k\}}\rangle$ $= \langle d_1+d_2-d_1d_2\rangle= \mathcal{R}_q
$; $S_{\{i_1,i_2,\cdots, i_k\}}\cap S'_{\{i_1,i_2,\cdots, i_k\}}$ $= \langle E_{\{i_1,i_2,\cdots, i_k\}}E'_{\{i_1,i_2,\cdots, i_k\}}\rangle= \langle 0\rangle$,
and  $S_{\{i_1,i_2,\cdots, i_k\}}+ S'_{\{i_1,i_2,\cdots, i_k\}}=\langle E_{\{i_1,i_2,\cdots, i_k\}}+E'_{\{i_1,i_2,\cdots, i_k\}}$ $-E_{\{i_1,i_2,\cdots, i_k\}}E'_{\{i_1,i_2,\cdots, i_k\}}\rangle= \langle e_1+e_2-e_1e_2\rangle=  \langle 1-\frac{1}{q}h \rangle$. This proves (i)-(iv).\vspace{2mm}

\noindent Using that $\frac{1}{q}h= (1-e_1-e_2)$ and $e_1e_2=0$ from Lemma 2 and noting that $e_1^2=e_1, e_2^2=e_2$ we find that $E_{\{i_1,i_2,\cdots, i_k\}} (\frac{1}{q}h)=0$.\\
 Similarly using $(e_1 +\frac{1}{q}h)=d_1$ and  $(e_2 +\frac{1}{q}h)=d_2$ from Lemma 2, we see that  $E_{\{i_1,i_2,\cdots, i_k\}}+ (\frac{1}{q}h)=D_{\{i_1,i_2,\cdots, i_k\}}$.\\ Therefore $S_{\{i_1,i_2,\cdots, i_k\}}\cap \langle \frac{1}{q}h\rangle = \langle E_{\{i_1,i_2,\cdots, i_k\}} (\frac{1}{q}h)
\rangle= \{0\},$ and $S_{\{i_1,i_2,\cdots, i_k\}}+ \langle \frac{1}{q}h\rangle  = \langle E_{\{i_1,i_2,\cdots, i_k\}} +\frac{1}{q}h-E_{\{i_1,i_2,\cdots, i_k\}}\frac{1}{q}h)\rangle=
\langle D_{\{i_1,i_2,\cdots, i_k\}}\rangle=Q_{\{i_1,i_2,\cdots, i_k\}}$. This proves (v) and (vi). \vspace{2mm}

\noindent Finally for $1 \leq k \leq [\frac{m}{2}] $, we have $$|Q_{\{i_1,i_2,\cdots, i_k\}}\cap Q'_{\{i_1,i_2,\cdots, i_k\}}|= | \langle\frac{1}{q}h\rangle|= p^m,$$ it being a repetition code over  $\mathcal{R}$. Therefore
$$ p^{mq}=|\mathcal{R}_q|=|Q_{\{i_1,i_2,\cdots, i_k\}}+ Q'_{\{i_1,i_2,\cdots, i_k\}}|=\frac{|Q_{\{i_1,i_2,\cdots, i_k\}}| |Q'_{\{i_1,i_2,\cdots, i_k\}}|}{|Q_{\{i_1,i_2,\cdots, i_k\}}\cap Q'_{\{i_1,i_2,\cdots, i_k\}}|}=\frac{|Q_{\{i_1,i_2,\cdots, i_k\}}|^2}{p^m}.$$
This gives $|Q_{\{i_1,i_2,\cdots, i_k\}}|=p^{\frac{m(q+1)}{2}}$.  Now   we find  that
$$ p^{\frac{m(q+1)}{2}}=|Q_{\{i_1,i_2,\cdots, i_k\}}|=|S_{\{i_1,i_2,\cdots, i_k\}}+ \langle\frac{1}{q}h\rangle|=|S_{\{i_1,i_2,\cdots, i_k\}}| |\langle\frac{1}{q}h\rangle|=|S_{\{i_1,i_2,\cdots, i_k\}}| p^m.$$
since $ |S_{\{i_1,i_2,\cdots, i_k\}}\cap \langle\frac{1}{q}h\rangle|= | \langle 0\rangle|= 1$. This gives $|S_{\{i_1,i_2,\cdots, i_k\}}|=p^{\frac{m(q-1)}{2}}$.  ~~~~~~~~~~~~~~~~~~~~~~~~$\Box$\vspace{2mm}

\noindent{\bf Theorem 5 :} If $p\equiv 1({\rm mod}~(m-1))$, $q\equiv  3({\rm mod}~4)$, $p$ is a quadratic residue modulo $q$, then
for all possible $\{i_1,i_2,\cdots, i_k\} \in \mathbb{A}$, the following assertions hold for quadratic residues codes over $\mathcal{R}$.
\vspace{2mm}\\
$\begin{array}{ll}
{\rm (i)} & Q_{\{i_1,i_2,\cdots, i_k\}}^{\perp} =  S_{\{i_1,i_2,\cdots, i_k\}},  \vspace{2mm}\\
{\rm (ii)}& S_{\{i_1,i_2,\cdots, i_k\}} {\rm ~are~ self~ orthogonal}. \end{array} $\vspace{2mm}

\noindent \textbf{Proof:} As $q\equiv  3({\rm mod}~4)$, $-1$ is a quadratic nonresidue modulo $q$. Therefore $j_1(x^{-1})=j_2$ and $j_2(x^{-1})=j_1$
and so $d_1(x^{-1})=1-e_1(x)$ and $d_2(x^{-1})=1-e_2(x)$. For $D_{\{i_1,i_2,\cdots, i_k\}}(x)=(\eta_{i_1}+\eta_{i_2}+\cdots +\eta_{i_k})d_1(x)+(1-\eta_{i_1}-\eta_{i_2}-\cdots -\eta_{i_k})d_2(x),$ $1- D_{\{i_1,i_2,\cdots, i_k\}}(x^{-1})=1-(\eta_{i_1}+\eta_{i_2}+\cdots +\eta_{i_k})(1-e_1(x))-(1-\eta_{i_1}-\eta_{i_2}-\cdots -\eta_{i_k})(1-e_2(x))= (\eta_{i_1}+\eta_{i_2}+\cdots +\eta_{i_k})e_1(x)+(1-\eta_{i_1}-\eta_{i_2}-\cdots -\eta_{i_k})e_2(x)=E_{\{i_1,i_2,\cdots, i_k\}}(x)$.  Now result (i) follows from Lemma 1. Using (vi) of Theorem 4, we have $S_{\{i_1,i_2,\cdots, i_k\}}\subseteq Q_{\{i_1,i_2,\cdots, i_k\}}= S_{\{i_1,i_2,\cdots, i_k\}}^{\perp}$. Therefore $S_{\{i_1,i_2,\cdots, i_k\}}$ are self orthogonal.  ~~~~~~~~~~~~~~~~~~~~~~~~$\Box$\vspace{2mm}

\noindent Similarly we get\vspace{2mm}

\noindent{\bf Theorem 6 :} If $p\equiv 1({\rm mod}~3)$, $q\equiv  1({\rm mod}~4)$, $p$ is a quadratic residue modulo $q$, then
for all possible choices of $\{i_1,i_2,\cdots, i_k\} \in \mathbb{A}$, the following assertions hold for quadratic residues codes over $\mathcal{R}$.\vspace{2mm}\\
$\begin{array}{ll}{\rm (i)} &Q_{\{i_1,i_2,\cdots, i_k \}}^{\perp} =  S'_{\{i_1,i_2,\cdots, i_k \}}, \vspace{2mm}\\ {\rm (ii)}& Q'^{\perp}_{\{i_1,i_2,\cdots, i_k\}} =  S_{\{i_1,i_2,\cdots, i_k\}}. \end{array} $\vspace{2mm}

 The extended QR-codes over $\mathbb{F}_{p}+u\mathbb{F}_{p}+\cdots+u^{m-1}\mathbb{F}_{p}$ are formed in the same way as the extended QR-codes over $\mathbb{F}_{p}$ are formed. \vspace{2mm}

\noindent{\bf Theorem 7 :} If  $q\equiv  3({\rm mod}~4)$, then
for all possible choices of $\{i_1,i_2,\cdots, i_k\} \in \mathbb{A}$, the extended QR-codes $\overline{Q_{\{i_1,i_2,\cdots, i_k\}}}$  of length $q+1$ are self dual.\vspace{2mm}

 \noindent \textbf{Proof:} If $p$ is odd, using quadratic reciprocity law, one easily finds that $-q$ is a quadratic residue modulo $p$ for  $q\equiv  3({\rm mod}~4)$. Then find an element $r\in \mathbb{F}_{p}$ such that $r^2\equiv -q({\rm~mod~} p)$. If $p=2$, $r^2\equiv -q({\rm~mod~} 2)$ clearly has a solution namely $r=1$. As $Q_{\{i_1,i_2,\cdots, i_k\}}=S_{\{i_1,i_2,\cdots, i_k\}}+\langle \frac{1}{q}h\rangle$, by Theorem 4, let $\overline{Q_{\{i_1,i_2,\cdots, i_k\}}}$ be the extended QR-code over $\mathcal{R}$ generated by
  $$~~~~~~~\begin{array}{cccccc}
     ~~~~~~~~~~~\infty & 0 & 1 & ~2 & \cdots & ~~q-1
   \end{array}\vspace{-2mm}$$ $$ \overline{G_{\{i_1,i_2,\cdots, i_k\}}}=\left(
                  \begin{array}{ccccc}
                    0 &  &  &  &  ~ \\
                    0 &  &  & G_{\{i_1,i_2,\cdots, i_k\}} &   \\
                    \vdots &  &  &  & ~  \\
                    r & ~1 & 1 &1 ~~\cdots ~~~~&  1
                  \end{array}
                \right)$$
 where $G_{\{i_1,i_2,\cdots, i_k\}}$ is a generator matrix for the even-like QR-code $S_{\{i_1,i_2,\cdots, i_k\}}$. The row above the matrix shows the column labeling  by $\mathbb{F}_q\cup \infty$.  Since the all one vector belongs to $Q_{\{i_1,i_2,\cdots, i_k\}}$ and its dual $Q_{\{i_1,i_2,\cdots, i_k\}}^{\perp}$ is equal to $  S_{\{i_1,i_2,\cdots, i_k\}}$, the last row of $\overline{G_{\{i_1,i_2,\cdots, i_k\}}}$ is orthogonal to all the previous rows of $\overline{G_{\{i_1,i_2,\cdots, i_k\}}}$. The last row is orthogonal to itself also as $r^2=-q $. Further as $S_{\{i_1,i_2,\cdots, i_k\}}$ is self orthogonal by Theorem 5, we find that the code $\overline{Q_{\{i_1,i_2,\cdots, i_k\}}}$ is self orthogonal. Now the result follows from the fact that $|\overline{Q_{\{i_1,i_2,\cdots, i_k\}}}|=p^m |S_{\{i_1,i_2,\cdots, i_k\}}|=p^{\frac{m(q+1)}{2}}= |\overline{Q_{\{i_1,i_2,\cdots, i_k\}}}^{\perp}|$. \vspace{2mm}

\noindent{\bf Theorem 8 :} If  $q\equiv 1({\rm mod}~4)$, then  for all possible choices of $\{i_1,i_2,\cdots, i_k\} \in \mathbb{A}$, the extended QR-codes satisfy $\overline{Q_{\{i_1,i_2,\cdots, i_k\}}}^{\perp}=\overline{Q'_{\{i_1,i_2,\cdots, i_k\}}}$.\vspace{2mm}

 \noindent \textbf{Proof:}  Let ~$\overline{Q_{\{i_1,i_2,\cdots, i_k\}}}$~ and ~$\overline{Q'_{\{i_1,i_2,\cdots, i_k\}}}$~ be the extended QR-codes over $\mathcal{R}$ generated by
  $$~~~~~~~\begin{array}{cccccc}
     ~~~~~~~~~~~~~~~~~~~\infty & 0 & 1 & ~2 & \cdots & ~~q-1
   \end{array}\vspace{-2mm}$$ $$ G_1=\overline{G_{\{i_1,i_2,\cdots, i_k\}}}=\left(
                  \begin{array}{ccccc}
                    0 &  &  &  &  ~ \\
                    0 &  &  & G_{\{i_1,i_2,\cdots, i_k\}} &   \\
                    \vdots &  &  &  & ~  \\
                    1 & ~1 & 1 &1 ~~\cdots ~~~~&  1
                  \end{array}
                \right)$$
 and  $$~~~~~~~\begin{array}{cccccc}
     ~~~~~~~~~~~~~~~~~~~\infty & 0 & 1 & ~2 & \cdots & ~~q-1
   \end{array}\vspace{-2mm}$$ $$G_2= \overline{G'_{\{i_1,i_2,\cdots, i_k\}}}=\left(
                  \begin{array}{ccccc}
                    0 &  &  &  &  ~ \\
                    0 &  &  & G'_{\{i_1,i_2,\cdots, i_k\}} &   \\
                    \vdots &  &  &  & ~  \\
                    -q & ~1 & 1 &1 ~~\cdots ~~~~&  1
                  \end{array}
                \right)$$ respectively where $ G_{\{i_1,i_2,\cdots, i_k\}}$ is a generator matrix for the  QR-code $S_{\{i_1,i_2,\cdots, i_k\}}$ and $ G'_{\{i_1,i_2,\cdots, i_k\}}$ is a generator matrix for the  QR-code $S'_{\{i_1,i_2,\cdots, i_k\}}$. Let $v$ denote the all one vector of length $q$. As $v \in Q'_{\{i_1,i_2,\cdots, i_k\}}$ and $Q'^{\perp}_{\{i_1,i_2,\cdots, i_k\}} =  S_{\{i_1,i_2,\cdots, i_k\}}$, $v$ is orthogonal to all the rows of $ G_{\{i_1,i_2,\cdots, i_k\}}$. Also $(-q,v)\cdot (1,v)=0$. Further rows of $ G'_{\{i_1,i_2,\cdots, i_k\}}$ are in $S'_{\{i_1,i_2,\cdots, i_k\}}=Q^{\perp}_{\{i_1,i_2,\cdots, i_k\}}$, so are orthogonal to rows of $ G_{\{i_1,i_2,\cdots, i_k\}}$. Therefore all rows of $G_2$ are orthogonal to all the rows of $G_1$. Hence  $\overline{Q'_{\{i_1,i_2,\cdots, i_k\}}} \subseteq \overline{Q_{\{i_1,i_2,\cdots, i_k\}}}^{\perp}$.  Now the result follows from comparing their orders.  \vspace{2mm}

\noindent \textbf{Corollary .} Let the matrix $V$ taken in the definition of the Gray map $\Phi$ satisfy $VV^T=\lambda I_m$, $\lambda \in \mathbb{F}_p^*$. If  $q\equiv  3({\rm mod}~4)$, then  for all possible choices of $\{i_1,i_2,\cdots, i_k\} \in \mathbb{A}$, the Gray images of extended QR-codes $\overline{Q_{\{i_1,i_2,\cdots, i_k\}}}$ ~ i.e.  $\Phi(\overline{Q_{\{i_1,i_2,\cdots, i_k\}}})$  are self dual codes of length $m(q+1)$ over $\mathbb{F}_p$ and the Gray images of the even-like QR-codes $S_{\{i_1,i_2,\cdots, i_k\}}$ i.e. $\Phi(S_{\{i_1,i_2,\cdots, i_k\}})$ are self-orthogonal codes of length $mq$ over $\mathbb{F}_p$. If  $q\equiv  1({\rm mod}~4)$, $\Phi(\overline{Q_{\{i_1,i_2,\cdots, i_k\}}})$  are formally self dual codes of length $m(q+1)$ over $\mathbb{F}_p$.\vspace{2mm}
\section{Examples}
The minimum distances of all the examples appearing in this section have been computed by the Magma Computational Algebra System.\\

\noindent {\bf Example 1 :} Take $m=3$, $p=7$, $q=3$ and $V =\left(
                                                               \begin{array}{ccc}
                                                                 2 & -2 & 1 \\
                                                                 1 & 2 & 2 \\
                                                                 2 & 1 & -2 \\
                                                               \end{array}
                                                             \right)
$ be a matrix over $\mathbb{F}_7$ satisfying $VV^T=2I$. The even like idempotent generators of quadratic codes of length 3 over $\mathbb{F}_7$
are $e_1=6x^2+3x+5$, $e_2=3x^2+6x+5$. The extended quadratic residue code $\overline{Q}_{\{1\}}$ is a self dual code over the ring $\mathbb{F}_7+u\mathbb{F}_7+u^2\mathbb{F}_7$ and its Gray image $\Phi(\overline{Q}_{\{1\}})$ is a self dual [12,6,4] code over $\mathbb{F}_7$. The Gray image of even like quadratic code  $S_{\{1\}}$ namely $\Phi(S_{\{1\}})$ is a self orthogonal nearly MDS [9,3,6] code over $\mathbb{F}_7$ in the sense of [2] as its dual having parameters [9,6,3] is also almost MDS. (Compare it with the example of Liu et al [6], where the code constructed was almost MDS but was not nearly MDS.) The Generator matrix of our nearly\vspace{2mm}\\ MDS code is $\left(
                                                      \begin{array}{ccccccccc}
                                                        1 & 0 & 0 & 0 & 4 & 5 & 6 & 3 & 2 \\
                                                        0 & 1 & 0 & 4 & 0 & 2 & 3 & 6 & 5 \\
                                                       0 & 0 & 1 & 5 & 2 & 3 & 2 & 5 & 3 \\
                                                      \end{array}
                                                    \right).$ It is better than the best known \vspace{2mm} \\linear code with these parameters (ref. Magma).\vspace{2mm}

\noindent {\bf Example 2 :} Take $m=3$, $p=5$, $q=11$ so that 5 is a quadratic residue modulo 11 and the matrix $V$ as in the above example.  The even like idempotent generators of quadratic codes of length 11 over $\mathbb{F}_5$
are $e_1=3x^{10}+x^9+3x^8+3x^7+3x^6+x^5+x^4+x^3+3x^2+x$, $e_2=x^{10}+3x^9+x^8+x^7+x^6+3x^5+3x^4+3x^3+x^2+3x$. The extended quadratic residue code $\overline{Q}_{\{1\}}$ is a self dual code over the ring $\mathbb{F}_5+u\mathbb{F}_5+u^2\mathbb{F}_5$ and its Gray image $\Phi(\overline{Q}_{\{1\}})$ is a self dual [36,18,9] code over $\mathbb{F}_5$. The Gray image of even like quadratic residue code  $S_{\{1\}}$  is a self orthogonal [33,15,10] code over $\mathbb{F}_5$.\vspace{2mm}

\noindent {\bf Example 3 :} Take $m=4$, $p=7$, $q=3$ and $$V =\left(
                                                               \begin{array}{cccc}
                                                                 2 & -2 & 1  & 1 \\
                                                                 -1 & 1 & 2&2 \\
                                                                 2 & 2 & 1&-1 \\1&1&-2&2
                                                               \end{array}
                                                             \right)
$$ be a matrix over $\mathbb{F}_{7}$ satisfying $VV^T=3I$. The even like idempotent generators of quadratic codes of length 3 over $\mathbb{F}_{7}$
are $e_1=6x^2+3x+5$, $e_2=3x^2+6x+5$.  The Gray image $\Phi(\overline{Q}_{\{1\}})$ of extended quadratic residue code $\overline{Q}_{\{1\}}$ is a  self dual [16,8,4] code over $\mathbb{F}_{7}$. The Gray image of even like quadratic residue code  $S_{\{1\}}$  is a self orthogonal [12,4,6] code over $\mathbb{F}_{7}$.\vspace{2mm}

\noindent {\bf Example 4 :} Take $m=5$, $p=13$, $q=3$ and $$V =\left(
                                                               \begin{array}{ccccc}
                                                                 2 & -2 & 1 &2 & -1 \\
                                                                 1 & 2 & 2&1&2 \\
                                                                 5 & 9 & -5&0&0 \\3&-5&-6&5&-6\\-1&2&0&7&-5
                                                               \end{array}
                                                             \right)
$$ be a matrix over $\mathbb{F}_{13}$ satisfying $VV^T=I$. The even like idempotent generators of quadratic codes of length 5 over $\mathbb{F}_{13}$
are $e_1=3x^2+x+9$, $e_2=x^2+3x+9$.  The Gray image $\Phi(\overline{Q}_{\{1\}})$ of extended quadratic residue code $\overline{Q}_{\{1\}}$ is a  self dual [20,10,4] code over $\mathbb{F}_{13}$. The Gray image of even like quadratic residue code  $S_{\{1\}}$  is a self orthogonal [15,5,6] code over $\mathbb{F}_{13}$.\vspace{2mm}

\noindent {\bf Example 5 :} Take $m=6$, $p=11$, $q=5$ and $$V =\left(
                                                               \begin{array}{cccccc}
                                                                 1 & 1 & 1 &1 & 1 & 1\\
                                                                 1 & 2 & -3&1&2&-3 \\
                                                                 1 & -3 & 2&1&-3&2 \\1&1&1&-1&-1&-1\\1&2&-3&-1&-2&3\\1&-3&2&-1&3&-2
                                                               \end{array}
                                                             \right)
$$ be a matrix over $\mathbb{F}_{11}$ satisfying $VV^T=6I$. The even like idempotent generators of quadratic codes of length 5 over $\mathbb{F}_{11}$
are $e_1=5x^4+8x^3+8x^2+5x+7$, $e_2=8x^4+5x^3+5x^2+8x+7$.  The Gray image $\Phi(\overline{Q}_{\{1\}})$ of extended quadratic residue code $\overline{Q}_{\{1\}}$ is a formally self dual [36,18,6] code over $\mathbb{F}_{11}$, whereas the Gray image $\Phi(S_{\{1\}})$ is a linear [30,12,8] code over $\mathbb{F}_{11}$.


\end{document}